\documentclass{jams-l}

\newtheorem{theorem}{Theorem}[section]

\theoremstyle{definition}

\theoremstyle{remark}

\numberwithin{equation}{section}
\usepackage{hyperref}
\usepackage{enumitem}

\makeatletter

\renewcommand\subsubsection{\@secnumfont}{\bfseries}%
\renewcommand\subsubsection{\@startsection{subsubsection}{3}
  \z@{.5\linespacing\@plus.7\linespacing}{-.5em}%
  {\normalfont\bfseries}}
  
  \makeatother

\begin{document}

\title{Guidelines for writing and reviewing software in computer algebra}


\author{Jeroen Hanselman}
\address{Gottlieb-Daimler Straße 48, 67655 Kaiserslautern}
\curraddr{}
\email{}
\thanks{}

\subjclass[2020]{Primary }

\date{}

\dedicatory{}

\begin{abstract}

The advent of computers has allowed mathematicians to do increasingly more difficult computations that used to be practically impossible. Peer reviewers will seldom look at any code attached to a math paper, however. In this article, we advocate for a software peer-reviewing process and will discuss some best practices for reviewers on how to peer review mathematical papers with code and providing guidelines for authors to improve the code that accompanies their mathematical paper.
\end{abstract}

\maketitle


\section{Why proper coding is important for mathematics}

The advent of computers has allowed mathematicians to do increasingly more difficult computations that used to be practically impossible. We now have databases of mathematical objects, measuring multiple terabytes in size, as well as huge computer algebra packages, proof assistants, and numerous other computational tools that aid us in doing mathematics. 

Despite the usefulness of computation in mathematics, it is not always used in a scientifically responsible way. Mathematicians have spent decades refining the way one should write down mathematics. Everything must be proven and written down clearly so that others can understand the logic that was used and agree that the results were correct. An article will most likely be rejected if it posits a theorem without proving it. Even though it is unlikely that a peer reviewer will study every detail and catch every mistake, they will try their best to see if the arguments and results look plausible and correct. When a paper relies on a computational part, however, the included code does not always get treated with the same care. On the one hand, reviewers might not always know how to properly review source code and, on the other hand, authors might not always know how to write good code. Because of this, the code might be badly structured, difficult to install and run, and sometimes even impossible to find. As a consequence, there are most assuredly papers floating around that contain proofs relying on source code that no one can run. This essentially means that these results come without a proof. This is not a problem that only occurs in mathematics, however. The replication crisis is a problem affecting many areas of science, is discussed, for example, in the article by Baker \cite{baker}. In computer science, people have started to do reviews for research artifacts, which are software or datasets that accompany a publication. It is discussed in articles such as \cite{Liu} and \cite{artifacts}, and we can see the reviewing of research artifacts being applied in practice at a conference like CHES, \url{https://ches.iacr.org/2024/artifacts.php}, where authors could have their code tested on Availability, (Re)usability and Reproducibility. 

As the mathematical objects we are dealing with become more complex, we are increasingly relying on computational software, and this trend will most likely continue. It is therefore crucial that we start thinking of guidelines and general practices that improve the reproducibility of our code. Otherwise, we might run into serious problems in the future.

Although we might prefer to have higher-quality code and a better peer-reviewing process for software, there are some good reasons why it has not happened yet. When doing research, software is mostly used as a tool. Just as we would use pen and paper. Although we are trained to transform the ideas we scribbled down somewhere into a rigorously written proof, we are not always taught to turn the code we wrote for ourselves into software that can be used and reused by the wider scientific community. A peer reviewer also does not always examine the code accompanying a paper. This could be due to a lack of time, lack of familiarity with the software used or a variety of other reasons. The code may also be in poor condition, making it a huge hassle for a reviewer to review. 

Furthermore, when applying for jobs, hiring committees do not always value the code you have in your repository as much as the number of publications you have. As the problem was recognized across many disciplines, it led to the creation of journals like JOSS (Journal of Open Source Software) \url{https://joss.theoj.org/}, which enables researchers to publish their software as a separate citable publication. 

Finally, writing good, reproducible software is hard and time-consuming. It is not just about writing legible code. One needs to think about a long-time storage solution, installation instructions, examples of how to use the software or data, and how to deal with different operating systems, among other things. The process is a lot more involved than just writing down a proof.

In summary, one sees that properly developing the software component of a paper is an important task, but that we, as a community, still need to do some work to integrate it into the scientific publication process. There are many ongoing initiatives to improve the way researchers deal with data. For example, the European Open Science Cloud (\url{https://eosc.eu/}) has a Task Force on "Infrastructure for Quality Research Software"; the United States Research Software Engineer Association (\url{https://us-rse.org/}) has a working group focused on Code Review; the Research Software Directory by the Helmholtz Centre (\url{https://helmholtz.software/}) promotes re-use and citation of research software; the EU's TIER2 project aims to improve reproducibility in science (\url{https://tier2-project.eu}) and FAIRsharing.org is a resource on data and metadata standards: \url{https://fairsharing.org}. The German government also recognized the need to adapt to the increasing reliance on computation in the sciences and established the NFDI (\url{https://www.nfdi.de}), an organization aimed at improving research data infrastructure across all sciences. I am part of the MaRDI project (short for Mathematical Research Data Initiative, \url{https://www.mardi4nfdi.de}), which is a subbranch of the NFDI specialized in mathematics. To read more about MaRDI's goals one can read the MaRDI white paper \cite{whitepaper}. As part of this project, I have been writing code reviews to try and figure out the best way to improve the quality of software written in mathematics. 

The following guidelines are based on my own experiences with reviewing code, and discussions with various authors and reviewers about best practices. I would like to thank John Cremona, Claus Fieker, Christiane Görgen, Lars Kastner, Jennifer Paulhus and Drew Sutherland for their helpful feedback. I am writing down what I have learned in the hopes of improving the community standards for handling code in mathematics. At this point, I have written approximately 100 software reviews for various conferences, most notably, for ANTS XVI (Algorithmic Number Theory Symposium, \url{https://antsmath.org/}) and LuCaNT (LMFDB, Computation, and Number Theory, \url{https://lucant.org/}) . Most of these reviews were in the area of algorithmic number theory. The introduction of a software reviewing process was done as an experiment, but the reactions have been positive. The inclusion of a software reviewing process seemed to have led to a significant increase in the amount of source code available compared to previous iterations of the ANTS conference. Authors were also happy about receiving feedback on the code that they wrote. A software reviewing process will most likely also be included in the next iterations of ANTS and LuCaNT.

\section{Some statistics about the impact of software reviewing}

To give an impression of the impact of software reviewing, I have compiled some statistics based on the reviews my MaRDI colleagues and I conducted for the ANTS and LuCaNT conferences. Do note that the conferences explicitly announced that software peer-reviewing would be part of the process, and authors were given the opportunity to submit their code if they had not done so initially. 

Henceforth, we will distinguish between papers and pieces of software. We will clarify the reason for this now. It is possible for a single paper to come with multiple pieces of software. One example of this is when a paper is accompanied by different implementations of the same algorithms, written in different programming languages. These usually came in separate repositories, which is another reason why they were reviewed separately. Some papers did not come with any code. Usually, when this happened, a link to the code was provided, but it was broken. If no link was given, and a paper clearly required an implementation to illustrate or prove the results discussed within, I treated it as if it came with a piece of software that I did not have access to. If the paper was of a more theoretical nature and the algorithms described did not really need an implementation to show their usefulness, I did not count the paper as being accompanied by a piece of software. Of the 75 pieces of software mentioned in the papers, 67 pieces of software  (89.3\%) had code available for review. In the other cases, including code with the paper was deemed essential, but the accompanying code was not available anywhere. Below, we will provide some statistics on the 67 pieces of software that came with code. To make the tables a bit more compact, we abbreviate the names of the categories: 
 \begin{itemize}
 \item \textbf{Lic}: This indicates that the software came with a license.
 \item \textbf{Rdme}: The code or repository had a Readme or a nice landing page explaining the contents of the repository.
\item \textbf{Repo}: The code was hosted on a repository for long-term storage.
\item \textbf{Inst(E)}: The code was easy to install (within 15 minutes).
\item \textbf{Inst(H)}: It was possible to install the code after putting in more effort.
\item \textbf{Specs}: The specifications necessary for allowing reproducibility were included in either the paper or the repository. (Specifications checked for included the version of software used, the CPU type, and the required amount of memory to run the code.) One common case where including specifications was deemed essential for reproducibility was when the algorithms were benchmarked and compared with other algorithms.
\item \textbf{Rdble}: The code in the software was judged to be ``readable'' in the sense of having proper indentation, well-chosen variable names and sufficient annotation, etc.
\item \textbf{Ran}: The code ran without any issues. (Regardless of whether the output actually matched with the claims made in the paper.)
\item \textbf{Correct}: The reviewer judged the results or data that were part of the publication were judged to be ``verifiable'' or ``likely correct''. 

\item \textbf{Ran and Correct}: If the code was judged to be \textbf{Correct} it did not necessarily mean that the code that produced the results ran flawlessly. Sometimes it was still possible to verify the results as correct in a different way. (When it was possible to circumvent a bug, for example, or if there was separate code that verified the results produced by the code that gave errors). On the flip side, just because the code ran, it does not necessarily mean that the output was correct. This is why we have a separate category for code that ran and also produced the correct results.

 \end{itemize}
\begin{table}[h]
\begin{tabular}{l|lllllll}
                   & Lic    & Rdme   & Repo   & Inst(E) & Inst(H) & Specs  & Rdble  \\
                   \hline
Pieces of Software & 29     & 48     & 64     & 52      & 59      & 34     & 49     \\
Percentage of Total\footnotemark[1] & 43.3\% & 71.6\% & 95.5\% & 77.6\%  & 88\%    & 50.7\% & 73.1\%
\end{tabular}
\end{table}
\footnotetext[1]{Total here means the 67 pieces of software for which code was available. } 
If we restrict ourselves to the 59 pieces of software that managed to install properly, we can look at how well the software ran:

\begin{table}[h]
\begin{tabular}{l|lll}
                   & Ran & Correct & Ran and Correct \\ \hline
Pieces of Software & 41       & 35                        & 30                                       \\
Percentage of Total\footnotemark[2] & 69.5\%   & 59.3\%                    & 50.8\%     
\end{tabular}
\end{table}
\footnotetext[2]{Total here means the 59 pieces of software for which the reviewer was able to install the code. } 

Although there is much room for improvement (only 30 out of 75 pieces of software had no problems at all), these numbers are honestly pretty good. If I were to increase the sample size and do software reviews for similar papers released 1-2 years ago, the numbers would most likely look much worse. Additionally, the problems we run into when trying to run the code are generally minor and can be fixed by an author within a reasonably short time. This in some way also highlights the advantages of a software reviewing process. It incentivizes authors to put more effort into their code, and it points out flaws in the authors' code that can be fixed before publication. 

To examine the effect the peer reviewing process had on the quality of the code, I went back and reviewed all of the papers that were accepted for publication again. I performed these reviews less thoroughly than the initial ones, primarily checking if the authors had made improvements based on feedback and if I was still able to get the code running.

 Of the 39 pieces of software mentioned in the papers that were accepted, about 36 pieces of software  (92.3\%) were submitted with the code used in the paper.

\begin{table}[h]
\begin{tabular}{l|lllllll}
                   & Lic    & Rdme   & Repo   & Inst(E) & Inst(H) & Specs  & Rdble  \\
                   \hline
Pieces of Software & 23     & 33     & 35     & 31      & 35      & 24     & 32     \\
Percentage of Total\footnotemark[3] & 63.8\% & 91.7\% & 97.2\% & 86.1\%  & 97.2\%    & 61.5\% & 82\%
\end{tabular}
\end{table}
\footnotetext[3]{Total here means the 36 pieces of software for which code was available. } 
If we restrict ourselves to the 35 pieces of software that managed to install properly we can look at how well the software ran:
\newpage

\begin{table}[h]
\begin{tabular}{l|lll}
                   & Ran & Correct & Ran and Correct \\ \hline
Pieces of Software & 32       & 29                        & 27                                       \\
Percentage of Total\footnotemark[4] & 91.4\%   & 82.8\%                    & 77.1\%     
\end{tabular}
\end{table}

Upon comparing the numbers, we find that the statistics of the accepted papers are much better than the submitted ones. We see that for 30 out of the 75 submitted pieces of software (40\%) there were no problems getting the software to run and verifying that the results were correct. Whereas this was true for 27 out of the 39 pieces of software (69.2\%) that were accepted for publication.

This marks a significant improvement. One possible reason for the improvement could be that having well-written software 
correlates with having a better written paper. Accepted papers tend to have higher-quality code, which leaves a better impression on peer reviewers and editors. Another reason is that many authors took the feedback from the code reviews to heart and actually made improvements to their code and their repositories. This ranged from small things, like adding a license or a short readme, to the addition of significant documentation, examples and the fixing of bugs. I saw positive changes for 19 out of the 35 pieces of software that were accepted. Overall, I believe the software reviewing process had a positive impact on the quality of the code attached to the accepted papers.

\section{Peer reviewing mathematics and peer reviewing research software.}

In this section, we will discuss which goals a software peer-reviewing process should aim to achieve. To start, we take a look at some aspects of traditional peer reviewing.

Greiffenhagen \cite{greiff} has written an article discussing the process and the role of peer reviewing in mathematics. We will shortly summarize some of his findings:
\begin{enumerate}
 \item Checking a mathematical proof for correctness is difficult and the peer-reviewing process cannot ensure that what is written down is 100\% correct. What it can do, however, is increase confidence in the correctness of the result. 

\item Peer reviewing a mathematical paper is generally more difficult and time-intensive than it is in other disciplines. As a result, the average number of reviewers that can look at a single paper is lower compared to those  disciplines.

\item The primary responsibility for ensuring that the stated results are correct lies with the author, but the reviewers share partial responsibility.

\item Reviewers generally do not read the entire article in detail. They skim the article to get a rough understanding of it and to determine whether they believe that the approach used could work, and then zoom in on the key arguments or parts that they do not understand to verify them for correctness.

\item Real mistakes are rarely found. It happens more often that reviewers ask for clarification for parts they do not understand.

\item It is a widely held belief that it is impossible to find all mistakes, but, if there are any, the community will eventually find them.

\end{enumerate}

Most of these points (with the possible exception of (2)) will likely also apply to any kind of software peer-reviewing process that we could design. It is impossible to ensure that the code will be completely bug-free (1), the authors of the code should bear the brunt of the responsibility to ensure that the code they published is correct (3), reviewers cannot be expected to check every line of code (4), it will be challenging for anyone to find mistakes in the code (5), and ideally most bugs will be discovered over time as others reuse the code, reimplement the algorithms, and so on (6).

In light of this the key goals of any peer-reviewing process should be the following:

\begin{enumerate}[label=\Alph*)]
\item The reviewing process should increase confidence that the results are correct.
\item The reviewing process should try to ensure that other mathematicians can follow the reasoning and reproduce the results by following the same steps as the authors.
\end{enumerate}

It is important to point out a few differences between having to review code and having to review a traditional mathematical paper.

First of all, the purpose of the code accompanying a paper may not necessarily be to prove a theorem. It is certainly common to outsource a part of a proof to a script written using a computer. It is then essentially used as a fancy calculator. But there are far more types of data and code that can come with a paper. A good way to think about it is to compare it to an experimental science like physics or chemistry. People do experiments to get data. These experiments and data are used as part of the publication in different ways:
\begin{itemize}
\item The experiments are used to construct some sort of database for future use (For example, the periodic table and all of the properties the elements contained within it have.)
\item The experiments are used to illustrate a part of the theory.
\item The data gathered from the experiments is used to formulate or support new theories.
\item The main goal of the paper is to explain the setup of a new kind of experiment.
\item The main goal of the paper is to explain how a certain kind of experiment can be improved to give faster/better results.
\end{itemize}

In computer algebra, we do very similar things. Only in our case, the experiments we conduct are performed using the algorithms that we write and use to study mathematics. To mirror the above list:
\begin{itemize}
\item We use algorithms to tabulate properties of mathematical objects. (See e.g. the atlas of finite groups \cite{atlas} or the LMFDB \cite{lmfdb}.)
\item We write code snippets to compute examples that are too difficult to compute by hand.
\item We perform computations to study conjectures and formulate new hypotheses.
\item We develop algorithms to compute things we were not able to compute before.
\item We write algorithms that improve on previous algorithms, making them faster and more efficient.
\end{itemize}

Consequently, to achieve goals A) and B), the peer reviewing process will need to differ slightly from the classical approach. To help a reviewer accomplish A) authors could, for example, add test files to the code that will help a reviewer check that there are no bugs in the code. Alternatively, the authors could add separate files that verify certain properties of the computed data. To help achieve B), authors must provide adequate installation instructions, hardware specifications, and documentation, ensuring that other mathematicians can reproduce their computations. 

Secondly, not all mathematicians will be able to understand code. Even if they are familiar with programming, they might not understand the specific programming language used in the program. Understanding a big repository filled with C code is far more difficult than executing a Sage script in a Jupyter Notebook, for example. This causes problems when trying to find a good reviewer. When looking for a reviewer of a traditional math paper, an editor tries to find an expert who understands the topic. If, in addition, it is also required that this expert understands the code and software used, the pool of potential reviewers will become extremely small. 

Because of this, we feel that when a paper with a significant software component has to be reviewed, an editor should split the peer review process into two tasks:
\begin{itemize}
\item  \textbf{Reviewing the mathematical content of the paper:}This should be done by someone who is an expert in the particular topic discussed by the paper.
\item \textbf{Reviewing the software that accompanied the paper:} This should be done by someone who is familiar with the programming languages involved. This person should also be broadly familiar with the field of research, but does not need to be an expert on the specific topic discussed in the paper. 
\end{itemize}

The reviewer and the code reviewer could potentially be the same person, and, if not, they could communicate if need be. This split will make it easier for editors to find reviewers as they do not necessarily need someone who is an expert in two things at the same time. As a consequence, the main focus of the code reviewer should not be to check for the correctness of the mathematics in the code. Their main focus should be on the reproducibility of the code. The average user should be able to find the software, install it, run it and check that the results are as the authors claim they are. A high level of reproducibility will ensure that other mathematicians will be able to use the code, check whether it is correct or not and figure out potential mistakes that were made. This approach has another advantage: it will reduce the time required to do a code review. And if it takes less time, it will make it easier to find willing reviewers.

We propose four categories that should be part of a software review:

\begin{enumerate}
  \item Metadata and the role of the code in the paper
  \item Reproducibility (Installation and rerunning the code)
  \item Reliability and Correctness
  \item Readability
\end{enumerate}

We will discuss these categories in more detail in the following sections, which will deal with mathematical software from the point of view of an author and from the point of view of a reviewer.

\section{General guidelines for writing reproducible code in Computer Algebra}

The next section will contain some guidelines and best practices for writing code for mathematical publications in computer algebra. The focus of an author seeking to make his research reproducible or a reviewer tasked with reviewing the code heavily depends on the role the code plays in the paper. We will briefly discuss some common scenarios before outlining the guidelines:

\begin{itemize}
    \item \textbf{The paper explains a new algorithm capable of computing things we were not able to compute before}: The main goal of the paper is to explain why this new algorithm works. The key part is to demonstrate its correctness. The theory behind the algorithm will primarily be explained in the paper, but examples and proper commentary in the repository will be useful for understanding the algorithm and verifying its correctness. Listing run times and system specifications will be less relevant in this case.
    \item \textbf{The paper presents an algorithm that is claimed to be better than existing algorithms:} In this case, providing hardware specifications and run times becomes crucial. If you claim your algorithm performs better than existing ones, you will need to provide benchmarking tools to support those claims.
\item \textbf{The paper discusses a database of mathematical objects:} It is difficult to provide proper documentation for a database of mathematical objects. Ideally, it should be possible to find out which people used which algorithms to produce each individual object in the database.
One would also like to know whether the data was computed rigorously or heuristically, and whether it was conditional on any conjectures or independent of them, as well as if the results have been double-checked using sanity checks or other reliability tests. Other things to consider are: How easy is it for users to access the database? Is there a user-friendly interface? Are there instructions for how to query the database? Is it possible for users to verify the data in the database themselves using some algorithms?

\item \textbf{The code is used to compute examples included in the paper:}
If the examples are not particularly hard to compute, just having a separate (clearly named) file for each example that immediately runs the needed code will suffice. If computing the examples will take a long time, it will be helpful to provide an estimate of the expected duration.

\item  \textbf{The code is used to do computations as part of a proof:} If the code is used as part of a proof, it will be crucial to clearly cross-reference between the code and the paper, as well as vice versa. Furthermore, the code should be properly commented, allowing readers of the proof to understand the software component. The code is an actual part of the proof and should be treated as such.
\end{itemize}

Now here are some general guidelines that may be helpful to keep in mind when writing research software.

\subsection{Metadata and the role of the code in the paper}

\subsubsection*{Make sure to describe the purpose of the accompanying code}
Write a few sentences in the paper clearly describing what kind of code accompanies the paper, what it was used for, and where it can be found. Based on my experience reviewing papers with software, I often found it difficult to find this information as it is usually scattered throughout the paper or missing altogether. Most people declare what the intent of their paper is in the abstract or the introduction, but it is also helpful to do something similar for the accompanying code.

\subsubsection*{Refer to sections of the paper in your code}
Conversely, it can be helpful to reference your paper in your code. This way, if someone does not understand what is happening, they can easily find the relevant sections in the paper.

\subsection{Reproducibility (Installation and rerunning the code)}

\subsubsection*{Make sure your code has a license}
It is important to ensure that the software you wrote comes with a license. Without a license, no one is legally allowed to reuse your code (which is detrimental to academia). If you have co-authors, they may technically be able to sue you if you reuse your own code, as they have partial ownership. For more information about what kind of license to use, check out: \href{https://choosealicense.com/}{https://choosealicense.com/} for example. For academia, it is usually best to pick some kind of open source license. Common ones are the GNU GPL license, the MIT license and the BSD licenses. One key difference is that if someone wants to reuse code that used a GPL license, the new work also needs to be open source and published under the same license. It can be used for commercial purposes, but everyone who receives a copy of the product must have access to the source code. The MIT and BSD licenses, however (which are relatively similar in practice) do not enforce that derivative works are published under the same license.  The main condition is that the new work contains proper attribution to the original author. The software may then also be turned into commercial closed-source software. If your code depends on software written by someone else, ensure the license you choose is compatible with theirs.

\subsubsection*{Make sure the code is publicly available in a long-lasting repository}
First, ensure the code is hosted on a platform where it can remain for multiple years without being deleted. A bad place to host would, for example, be on your university's or your own website. It happens often that when people change jobs or leave academia the corresponding websites disappear. Even if you move the code elsewhere, papers that cite the code will still reference the old location, making it difficult to find the code. It is therefore recommended to host your code on a public platform like Github, Gitlab or Zenodo. Github has the advantages of having a lot of features and being the most popular platform for hosting code. One possible negative is that Github is a subsidiary of Microsoft. However unlikely, they might choose to monetize Github or impose some other unwanted changes at some point. Zenodo is hosted by the CERN data center and more specifically geared towards science. Every upload also has a DOI (a persistent identifier), which makes them more easily citable, just like papers.

\subsubsection*{Ensure your code has a Readme file}

The goal of the Readme file should be to serve as the visiting card of your repository. It should give readers a quick impression of what the purpose of the repository is and what one can do with it. It is furthermore a good place to provide an overview of the files contained in the repository and their purpose. The Readme could also potentially include installation instructions, basic documentation,  and some small examples. Finally, it can be used to reference the accompanying paper and cite works by other authors that your software depends on.

\subsubsection*{Provide clear installation instructions}

Other mathematicians need to be able to run your code to verify your results, gain a deeper understanding of your work, or reuse it and improve upon it. Unfortunately, many repositories currently do not include any instructions at all or only very vague ones. The average user will not always be capable of finding out exactly what they need to do, so make the instructions as explicit as possible. A good test would be to let a colleague or a friend try to install the software you wrote and see if they manage. Sometimes the issue is not with your own code, but with the installation instructions someone else wrote (or did not write) for a package you depend on. In that case, you should either try to contact the original authors of that code and ask them to improve their instructions or provide better instructions on your own page.

\subsubsection*{Think about what you could do to make the installation as easy as possible}

This is strongly related to the previous point. While complicated installation instructions are sometimes unavoidable, there are steps that can be taken to make the process easier for users. It may be possible to write an installation script that automatizes a couple of steps. Additionally, try to minimize the need for users to manually change files as much as possible.

Unfortunately, I have seen  it happen multiple times where the installation instructions require users to install a library from a different author and then manually modify certain lines in that package. This is a bad practice that can cause a lot of issues. Will the library still work with your changes? What if the library gets updated in the future? 

\subsubsection*{Provide documentation for the code}
Other users should be able to understand how to use your code. Every function you write should include an explanation of what the input and output, as well as a short description of its purpose. Also make sure you include a few files with examples of how to use your code. Ideally, your code should include enough explanatory text to help users understand what it does. A Jupyter notebook would also work well for something like this.

\subsubsection*{Mention the specifications of the system you used to do your computations}
Be sure to list key hardware specifications (CPU, RAM), the OS used and the exact versions of the programming languages and packages involved. For large or difficult computations, it may also be helpful to record the amount of data produced and the time it took for the computations to finish.

\subsubsection*{Citation}
Make sure to give proper credit to the work you are building on. Cite the software and the packages written by other authors that you are using in your repository and in the accompanying paper.

\subsection{Correctness and reliability}

\subsubsection*{Tests} To boost confidence in the correctness of your code, it is good practice to write some test code. The test code does not necessarily need to test your intended example. Its goal could be to test certain functions or submodules of the entire package, perform sanity checks, verify the code is working in certain edge cases or to test randomly generated examples.

Running all the code on a real example could take a long time, but cleverly chosen tests can make a convincing case that the code works in principle despite taking less time to run. When writing test code, aim to cover as much code as possible. There also exist tools measuring code coverage that can do this for you. 

Test code can also help you prevent unwanted bugs from sneaking in. If you change some code, and your tests no longer work as a result, you know that you have messed something up. To help with this, tools for continuous integration also exist. These tools can help you run tests automatically every time you change something in the code.

\subsubsection*{Verification files} If you want to prove to someone that your computations are correct, you might not always want to redo the computation, as it could take a lot of time. However, verifying that something is correct isn't necessarily similar to redoing the computation and can sometimes be done in a different and potentially quicker way. Computing the inverse of a matrix $A$, for example, is different from multiplying $A$ with the matrix you computed and checking if the result is the identity matrix. These verification tests also increase confidence that your algorithm is correct.

\subsubsection*{Multiple implementations} To check whether the data you computed is correct, it can be helpful to implement your code in multiple software systems. In some cases, different systems actually produce different results (due to bugs in the software systems). Testing the implementation in different systems increases the chances of your results being correct.

\subsection{Readability}

\subsubsection*{Use clear naming conventions}
It might be a small thing, but thinking a bit longer about proper names for files, functions and variables will do wonders for making it more comprehensible for other people (and even for yourself).  No one likes to see names like \textit{helperfunctions3} or \textit{mmhm}. 
\begin{itemize}
    \item Try to make your naming conventions meaningful, so people who have read the paper can immediately understand what each term represents.
    \item Try to stick to common mathematical variables such as $f$ for function and $v$ for vector, etc.
    \item Try to use as few abbreviations as possible. Longer names might help with legibility even if it takes a bit longer to type.
    \item Try to remain consistent with your naming scheme.
    
\end{itemize}

\subsubsection*{Ensure your repository is clearly structured and that files are easy to find} This is especially important if the repository is big. If there are hundreds of files, you don't want people to spend hours trying to figure out what each one does and where to find what they need.

\section{General guidelines for reviewing code in Computer Algebra}

Here we will discuss an approach to reviewing code for a mathematical paper. The first thing to do is to read the abstract and the introduction to get a rough impression of what the article is about. After this, it makes sense to skim through the rest of the article and to take a look at the files provided. This should give you sufficient information to determine all the necessary details for the first part of the review, which we will discuss in the following section.

\subsection{Metadata and the role of the code in the paper}

The review should contain the following metadata and information about the paper and its accompanying code:
\begin{itemize}
 \item The title of the paper.
 \item The authors' names (if given).
 \item The date on which the review was written.
 \item An overview of the files accompanying the paper, such as: code, notebooks, computed data, example files, documentation, Docker files/virtual machines, and test files.
 \item The version of the code that was reviewed. (If version information is not available, the download date of the code can be listed.)
 \item The link that was used to download the code. 
 \item A short description of the relationship between the code and the paper. What does the code add to the paper? How important is it for the results? 
\end{itemize}

After the first phase of the software review, you should have a good idea of what to expect to find in the code. Ideally, it should be stored somewhere in a repository (or multiple repositories.) So the next step would be to take a look at the repository and its structure.

\subsection{Reproducibility}

A quick scan of the repository can already reveal whether information essential for ensuring reproducibility is included. We give a rough checklist of things to look at:

\begin{itemize}
  \item Does the repository have a license?
  \item Where is the code hosted? The author's own personal website, for example, is a bad idea as it is not guaranteed to be there for a long time. A platform like Zenodo would be ideal, but a Github repository is still fine.
  \item Is the code publicly available? 
  \item Is the code open source? 
  \item Does the repository have a DOI that links to it? 
  \item Does the repository have a Readme? 
  \item What is the quality of this Readme? 
\end{itemize}

At this point, the reviewer should attempt to install the code and run some of the examples or test code included in the repository or paper. Write down all the steps you've taken when trying to install the code. If the installation fails,  the author will be able to use this information to help identify what caused the error.

Some things that could be included in a review:

\begin{itemize}
  \item Give a general impression of how difficult it was to install the software. Is it doable for the average user? Could it be improved? 
  \item Can the installation instructions be improved?
  \item A list of all the steps the reviewer took to try to get the code running.
  \item Check what information the authors gave about the setup they used. This includes CPU, RAM, what operating system was used, what software was used and version information. It can also include bench-marking information and running times used by the authors.
  \item See if the authors properly cited all the software they use in their program.
\end{itemize}

To help the author process the feedback he gets from the review, it is a good idea to also provide some information of the setup used for the review. Be sure to list:
\begin{itemize}
\item What operating system the reviewer used for the review.
 \item What programming languages or standard software packages were used for the review. (Including version information).
\end{itemize}

Once the software is installed and runs successfully on some basic examples, the reviewer should attempt to reproduce the authors' work as described in their paper.

\subsection{Reliability and correctness}
There are various ways to increase confidence in the correctness of the results stated in the paper. The main way for the reviewer to test this is to try to run the code itself and see if the results match the authors' claims. However, this might not always be feasible. For instance, if the computations the authors did took months, it is unrealistic to expect the reviewer to rerun them. Nevertheless, there are some solutions that can increase confidence that running the code on these big examples will also produce the correct results. 

\begin{itemize}
\item The repository could include smaller examples and show the code works correctly there.
\item The repository could contain test files that show that parts of the code are working correctly.
\item Sometimes verifying that the results are correct is easier and less time-intensive than doing the computations. The repository could include files that verify the correctness of, or perform some sanity checks on, the computed results.
\end{itemize}

In any case, if an error arises or the results are different from what was expected, the reviewer should describe exactly what they did and what went wrong. Sometimes it is easy to identify the issue, and the reviewer could attempt a quick fix, but in principle, it is not the reviewer's task to debug the authors' code.

\subsection{Readability}

The final thing a reviewer should check is whether the code was written in such a way that someone who would really want to study the code in detail for debugging or improving it could feasibly do so. This does not mean that the reviewer should try to understand the code. But they should check if some basic good coding practices are being respected.
\begin{itemize}
\item Do functions have a description of what they do? Is it clearly described what the input and output of the function are?
\item Is the code properly annotated?
\item Does the indentation look okay?
\item Are files, functions and variables given distinctive and descriptive names?
\item Is the repository structured sensibly?
\end{itemize}

A software review that takes the four categories discussed above into account should provide a fairly comprehensive overview of the quality of the code without being too time-consuming for the reviewer doing the review. 

An example of the reviewing template I have used for reviewing code can be found here: \url{https://github.com/JHanselman/code-review-template}

\section{Future outlook}

I personally consider it extremely important for us as a mathematical community to start implementing some form of quality control for research software and research data in mathematics. There will not be a solution that fits everyone, however. What might be important for a computer algebraist might not be as important to an applied mathematician, for example. These guidelines are a step towards getting people to think and discuss about what kind of solutions would work best. 

And discussion is what we need. I hope anyone who thinks this is important will talk to their colleagues, talk to their students, talk to editors of journals, organizers of conferences, and whoever can make a change in how we deal with research data and research software. If we want to improve the situation, achieve higher-quality reusable code and data, and gain recognition for our coding efforts, the change must come from within our community. We need to decide that we want this, and we need to actively start working towards it.

\bibliographystyle{amsplain}
\bibliography{references.bib}

\end{document}